\documentclass[12pt, reqno]{amsart}
\usepackage{amsmath, amsthm, amscd, amsfonts, amssymb, graphicx}
\usepackage[bookmarksnumbered,plainpages]{hyperref}
\input{mathrsfs.sty}

\newtheorem{theorem}{Theorem}[section]

\newtheorem{proposition}[theorem]{Proposition}
\newtheorem{corollary}[theorem]{Corollary}
\theoremstyle{definition}

\theoremstyle{remark}
\newtheorem{remark}[theorem]{Remark}

\numberwithin{equation}{section}

\begin{document}
\title{An operator extension of Bohr's inequality}
\author[M.S. Moslehian, J. Pe\v{c}ari\'{c}, I. Peri\'c]{Mohammad Sal Moslehian$^1$, Josip Pe\v{c}ari\'{c}$^2$ and
Ivan Peri\'c$^3$}
\address{$^1$ Department of Pure Mathematics, Ferdowsi University of Mashhad, P.O. Box 1159, Mashhad 91775, Iran;
\newline Center of Excellence in Analysis on Algebraic Structures
(CEAAS), Ferdowsi University of Mashhad, Iran.}
\email{moslehian@ferdowsi.um.ac.ir and moslehian@ams.org}
\address{$^2$ Faculty of Textile Technology, University of Zagreb, Pierottijeva 6, 10000 Zagreb, Croatia }
\email{pecaric@mahazu.hazu.hr}
\address{$^3$ Faculty of Food Technology and Biotechnology, Mathematics
Department, University of Zagreb, Pierottijeva 6, 10000 Zagreb,
Croatia} \email{iperic@pbf.hr} 


\begin{abstract} We establish an operator extension of the following generalization of Bohr's
inequality, due to M.P.~Vasi\'c and D.J.~Ke\v{c}ki\'{c}:
$$\left|\sum_{i=1}^n z_i\right|^r \leq \left(\sum_{i=1}^n
\alpha_i^{1/(1-r)}\right)^{r-1}\sum_{i=1}^n \alpha_i|z_i|^r \quad (
r>1, z_i \in{\mathbb C}, \alpha_i>0, 1 \leq i \leq n)\,.$$ We also
present some norm inequalities related to our noncommutative
generalization of Bohr's inequality.

\end{abstract}
\maketitle


\section{Introduction}

Let ${\mathfrak A}$ be a $C^*$-algebra of Hilbert space operators
and let $T$ be a locally compact Hausdorff space. A field
$(A_t)_{t\in T}$ of operators in ${\mathfrak A}$ is called a
continuous field of operators if the function $t \mapsto A_t$ is
norm continuous on $T$. If $\mu$ is a Radon measure on $T$ and the
function $t \mapsto \|A_t\|$ is integrable, one can form the Bochner
integral $\int_{T}A_t{\rm d}\mu(t)$, which is the unique element in
${\mathfrak A}$ such that
$$\varphi\left(\int_TA_t{\rm d}\mu(t)\right)=\int_T\varphi(A_t){\rm d}\mu(t)$$
for every linear functional $\varphi$ in the norm dual ${\mathfrak
A}^*$ of ${\mathfrak A}$; cf. \cite[Section 4.1]{H-P}.\\
Further, a field $(\varphi_t)_{t\in T}$ of positive linear mappings
$\varphi: {\mathfrak A} \to {\mathfrak B}$ between $C^*$-algebras of
operators is called continuous if the function $t \mapsto
\varphi_t(A)$ is continuous for every $A \in {\mathfrak A}$. If the
$C^*$-algebras include the identity operators, denoted by the same
$I$, and the field $t \mapsto \varphi_t(I)$ is integrable with
integral $I$, we say that $(\varphi_t)_{t \in T}$ is unital.

The classical Bohr's inequality states that for any $z, w \in
{\mathbb C}$ and any positive real numbers $r, s$ with
$\frac{1}{r}+\frac{1}{s}=1$,
$$|z+w|^2 \leq r |z|^2 + s
|w|^2.$$ This inequality admits the operator extension
$$|A+B|^2 \leq r|A|^2+s|B|^2$$
for operators $A, B$ in the algebra ${\mathbb B}({\mathscr H})$ of
all bounded linear operators on a complex Hilbert space ${\mathscr
H}$ (to see this use the Cauchy--Schwarz inequality and the fact
that the operator $C$ is positive if and only if $\langle
Cx,x\rangle\geq 0$).

\noindent Over the years, interesting generalizations of this
inequality have been obtained in various settings; cf. \cite{C-P,
HIR, M-R, P-S, P-R, ZHA}. There is one of special interest given by
M.P.~Vasi\'{c} and D.J.~Ke\v{c}ki\'{c} \cite{V-K1}:

\noindent If $z_1, \cdots, z_n$ are complex numbers, $r>1$ and
$\alpha_i >0 \, (i=1, 2, \cdots, n)$, then
\begin{eqnarray}\label{class}
\left|\sum_{i=1}^n z_i\right|^r \leq \left(\sum_{i=1}^n
\alpha_i^{1/(1-r)}\right)^{r-1}\sum_{i=1}^n \alpha_i|z_i|^r\,.
\end{eqnarray}
This is indeed an immediate consequence of the H\"older inequality.
In this paper we establish the operator version of inequality
\eqref{class} and apply the obtained operator inequalities to obtain
some norm inequalities related to our operator extension of Bohr's
inequality.


\section{Main results}

Recall that a continuous real function $f$ defined on a real
interval $J$ of any type is said to be operator convex if $f(\lambda
A+(1-\lambda)B) \leq \lambda f(A)+ (1-\lambda)f(B)$ holds for all
$\lambda \in [0,1]$ and all self-adjoint operators $A, B$ acting on
a Hilbert space with spectra in $J$. For instance $f(x)=x^r$, where
$1 \leq r \leq 2$ is operator convex on $[0,\infty)$; see \cite[p.
123]{BHA}. In \cite{H-P-P}, the authors gave a general formulation
of Jensen's inequality for unital fields of positive linear mappings
in which they dealt with operator convex functions.

\noindent We need the main result \cite[Theorem 2.1]{H-P-P}. We
state it for the sake of convenience.
\begin{theorem}\label{th-1}
Let $f$ be an operator convex function on an interval $J$ and let
${\mathfrak A}$ and ${\mathfrak B}$ be unital $C^*$-algebras. If
$(\varphi_t)_{t\in T}$ is a unital field of positive linear mappings
$\varphi_t: {\mathfrak A} \to {\mathfrak B}$ defined on a locally
compact Hausdorff space $T$ with a bounded Radon measure $\mu$, then
the inequality
\begin{eqnarray*}
f\left(\int_T \varphi_t(A_t)d\mu(t)\right) \leq
\int_T\varphi_t(f(A_t))d\mu(t)\,.
\end{eqnarray*}
holds for every bounded continuous field $(A_t)_{t\in T}$ of
self-adjoint elements of ${\mathfrak A}$ with spectra contained in
$J$.
\end{theorem}
Utilizing the theorem above we prove our main result.
\begin{theorem}\label{Th1}
Let ${\mathfrak A}$ and ${\mathfrak B}$ be $C^*$-algebras of
operators containing $I$, let $T$ be a locally compact Hausdorff
space equipped with a bounded Radon measure $\mu$, let
$\left(\alpha_t\right)$ a bounded continuous nonnegative function
such that $\left(\alpha_t\right) \in L^{\frac{1}{1-r}}(T,\mu)$ and
$1 < r \leq 2$. Let $(A_t)_{t\in T}$ be a bounded continuous field
of positive elements in ${\mathfrak A}$ and let $(\varphi_t)_{t\in
T}$ be a field of positive linear mappings $\varphi_t: {\mathfrak A}
\to {\mathfrak B}$ defined on $T$ satisfying
\begin{eqnarray*}
\int_T\alpha_t^{1/(1-r)}\varphi_t(I)d\mu(t) \leq
\int_T\alpha_t^{1/(1-r)}d\mu(t) I\,.
\end{eqnarray*}
Then
\begin{eqnarray}\label{Pec}
\left(\int_T \varphi_t(A_t)d\mu(t) \right)^r \leq \left(\int_T
\alpha_t^{1/(1-r)}d\mu(t)\right)^{r-1}\int_T
\alpha_t\varphi_t(A_t^r)d\mu(t)\,.
\end{eqnarray}
\end{theorem}
\begin{proof}
Let $\infty$ be an object not belonging to $T$. Consider
$T_\infty:=T\cup\{\infty\}$ as a locally compact topological space
by equipping $\{\infty\}$ with the discrete topology and extend
$\mu$ on $T_\infty$ by $\mu(\{\infty\})=1$. Set $f(x)=x^r\,\, (x\in
[0,\infty))$, $\tilde{\varphi}_t:=\frac{P_t}{Q} \varphi_t\,\, (t\in
T_{\infty})$, where $P_t:=\alpha_t^{1/(1-r)}$, $P_{\infty}:=1$,
$Q:=\int_T\alpha_t^{1/(1-r)}d\mu(t)$ and
$$\varphi_{\infty}(A):=\langle
Ae,e\rangle\left(\int_TP_td\mu(t) I -
\int_TP_t\varphi_t(I)d\mu(t)\right)\,,$$ in which $A$ belongs to the
$C^*$-algebra ${\mathfrak A}$ acting on a Hilbert space ${\mathscr
H}$ and $e \in {\mathscr H}$ is a fixed unit vector. Then
\begin{eqnarray*}
\int_{T_\infty}\tilde{\varphi}_t(I)d\mu(t)&=&\frac{1}{Q}\int_TP_t\varphi_t(I)d\mu(t)+
\tilde{\varphi}_{\infty}(I)\\
&=&\frac{1}{Q}\int_T P_t\varphi_t(I)d\mu(t)+\frac{1}{Q}\left(\int_T P_td\mu(t)I-\int_T P_t\varphi_t(I)d\mu(t)\right)\\
&=&I\,.
\end{eqnarray*}
It follows from Theorem \ref{th-1} that
\begin{eqnarray*}
\left(\int_{T_{\infty}}
\tilde{\varphi}_t(\tilde{A}_t)d\mu(t)\right)^r \leq
\int_{T_{\infty}}\tilde{\varphi}_t(\tilde{A}_t^r)d\mu(t)
\end{eqnarray*}
namely
\begin{eqnarray}\label{Pec2}
\left(\int_{T_{\infty}} P_t\varphi_t(\tilde{A}_t)d\mu(t)\right)^r
\leq Q^{r-1} \int_{T_{\infty}}P_t\varphi_t(\tilde{A}_t^r)d\mu(t)
\end{eqnarray}
for all bounded continuous fields $(\tilde{A_t})_{t\in T_\infty}$.\\
Put $\tilde{A}_t=A_t/P_t\,\, ( t\in T)$ and $\tilde{A}_{\infty}=0$
in (\ref{Pec2}) to obtain
\begin{eqnarray}\label{Pec3}
\left(\int_T \varphi_t(A_t)d\mu(t)\right)^r \leq
Q^{r-1}\int_TP_t^{1-r}\varphi_t(A_t^r)d\mu(t)\,.
\end{eqnarray}
By the definitions of $P_t$ and $Q$, \eqref{Pec3} becomes
\eqref{Pec}.
\end{proof}

\begin{remark} Using analogous arguing as in the proof of Theorem
\ref{Th1}, one can prove the following Jensen's inequality. Suppose
that $\left(\varphi_t\right)$ is a continuous field of positive
linear mappings $\varphi_t: {\mathfrak A}\to {\mathfrak B}$,
${\mathfrak A}$ and ${\mathfrak B}$ are unital $C^*$-algebras,
$\left(A_t\right)$ is a bounded continuous field of self-adjoint
elements in ${\mathfrak A}$ with spectra in an interval $J$ such
that $0\in J$, $(\beta_t)$ is a continuous nonnegative function such
that $\int_{T}\beta_t d\mu (t)>0$ and
\begin{equation}\label{GenOpJensen0}
\int_{T}\beta_t\varphi_t\left(I\right)d\mu (t)\leq \int_{T}\beta_t
d\mu (t)I\,.
\end{equation}
If  $f$ is an operator convex function on an interval $J$ such that
$f(0)\leq 0$, then
\begin{equation}\label{GenOpJensen1}
f\left(\frac{1}{\int_{T}\beta_td\mu(t)}\int_T\beta_t\varphi_t\left(A_t\right)d\mu
(t)\right)\leq \frac{1}{\int_{T}\beta_td\mu
(t)}\int_T\beta_t\varphi_t\left(f\left(A_t\right)\right)d\mu
(t).\end{equation} If equality holds in (\ref{GenOpJensen0}), then
it is not necessary to assume that $0\in J$ and $f(0)\leq 0$.
\end{remark}
A discrete version of the theorem above is the following result
obtained by taking $T=\{1, \cdots, n\}$.

\begin{corollary}\label{new}
Let $1 < r \leq 2$, $\alpha_i >0\,\,(i=1, \cdots, n)$, let $A_1,
\cdots A_n$ be positive operators acting on a Hilbert space
${\mathscr H}$ and let $\varphi_i\,\, (i=1, \cdots, n)$ be positive
linear mappings on ${\mathbb B}({\mathscr H})$ satisfying
\begin{eqnarray}\label{ppeecc}
\sum_{i=1}^{n}\alpha_i^{1/(1-r)}\varphi_i(I) \leq
\sum_{i=1}^{n}\alpha_i^{1/(1-r)} I\,.
\end{eqnarray}
Then
\begin{eqnarray*}
\left(\sum_{i=1}^n \varphi_i(A_i) \right)^r \leq \left(\sum_{i=1}^n
\alpha_i^{1/(1-r)}\right)^{r-1}\sum_{i=1}^n
\alpha_i\varphi_i(A_i^r)\,.
\end{eqnarray*}
\end{corollary}

\noindent By setting $\varphi_i(A)=X_i^*AX_i$ in Corollary \ref{new}
we find that
\begin{corollary}
Let $1< r \leq 2$, $\alpha_i >0$ for $i=1, \cdots, n$, let $A_1,
\cdots A_n$ be bounded operators acting on a Hilbert space
${\mathscr H}$ with $A_i \geq 0$ and let $X_1, \cdots, X_n \in
{\mathbb B}({\mathscr H})$ satisfying
\begin{eqnarray}\label{cond}
\sum_{i=1}^{n}\alpha_i^{1/(1-r)}X_i^*X_i \leq
\sum_{i=1}^{n}\alpha_i^{1/(1-r)} I\,.
\end{eqnarray}
Then
\begin{eqnarray*}
\left(\sum_{i=1}^n X_i^*A_iX_i\right)^r \leq \left(\sum_{i=1}^n
\alpha_i^{1/(1-r)}\right)^{r-1}\sum_{i=1}^n \alpha_iX_i^*A_i^rX_i\,.
\end{eqnarray*}
\end{corollary}


\noindent Condition (\ref{cond}) trivially holds if $X_i^*X_i \leq
I$ for all $i=1, \cdots, n$. In fact we can give another proof of
the result in this case.

\begin{corollary}\label{xi}
Let $A_1, \cdots A_n, X_1, \cdots, X_n \in {\mathbb B}({\mathscr
H})$ with $A_i \geq 0$, $X_i^*X_i \leq I$ for $i=1, \cdots, n$, and
let $1< r \leq 2$, $\alpha_i >0$ for $i=1, \cdots, n$. Then
\begin{eqnarray*}
\left(\sum_{i=1}^n X_i^*A_iX_i\right)^r \leq \left(\sum_{i=1}^n
\alpha_i^{1/(1-r)}\right)^{r-1}\sum_{i=1}^n \alpha_iX_i^*A_i^rX_i\,.
\end{eqnarray*}
\end{corollary}
\begin{proof} First note that $0 \leq \left(X_i^*A_iX_i\right)^r \leq X_i^*A_i^rX_i$
for each $i$; cf. \cite[Theorem 2.1]{HAN}. For $i=1, \cdots, n$, set
$\beta_i=\alpha_i^{1/(1-r)}$, $B_i=X_i^*A_iX_i/\beta_i$. We have
\begin{eqnarray*}
\left(\sum_{i=1}^n X_i^*A_iX_i\right)^r &=& \left(\sum_{i=1}^n \beta_iB_i\right)^r\\
&=&\left(\sum_{i=1}^n \beta_i\sum_{j=1}^n \frac{\beta_j}{\sum_{k=1}^n \beta_k}B_j\right)^r\\
&=& \left(\sum_{i=1}^n \beta_i\right)^r\left(\sum_{i=1}^n \frac{\beta_i}{\sum_{k=1}^n \beta_k}B_i\right)^r\\
&\leq& \left(\sum_{i=1}^n \beta_i\right)^r\frac{\sum_{i=1}^n \beta_iB_i^r}{\sum_{i=1}^n \beta_i}\\
&&\qquad(\textrm{by the operator convexity of $f(t)=t^r$})\\
&=& \left(\sum_{i=1}^n \beta_i\right)^{r-1}\sum_{i=1}^n \beta_i^{1-r}(X_i^*A_iX_i)^r\\
&=& \left(\sum_{i=1}^n \alpha_i^{1/(1-r)}\right)^{r-1}\sum_{i=1}^n \alpha_i (X_i^*A_iX_i)^r\\
&\leq& \left(\sum_{i=1}^n
\alpha_i^{1/(1-r)}\right)^{r-1}\sum_{i=1}^n \alpha_i
X_i^*A_i^rX_i\,.
\end{eqnarray*}
\end{proof}


\begin{proposition}\label{2}
Let $A_1, \cdots A_n \in {\mathbb B}({\mathscr H})$ with
$A_i^*A_j=0$ for $1 \leq i \neq j \leq n$, and let $2< r \leq 4$,
$\alpha_i >0$ for $i=1, \cdots, n$. Then
\begin{eqnarray*}
\left|\sum_{i=1}^n A_i\right|^r \leq \left(\sum_{i=1}^n
\alpha_i^{2/(2-r)}\right)^{(r-2)/2}\sum_{i=1}^n \alpha_i|A_i|^r
\end{eqnarray*}
and
\begin{eqnarray}\label{norm}
\left\|\sum_{i=1}^n A_i \right\|^r \leq \left(\sum_{i=1}^n
\alpha_i^{2/(2-r)}\right)^{(r-2)/2}\sum_{i=1}^n \alpha_i\|A_i\|^r.
\end{eqnarray}
\end{proposition}
\begin{proof}
\begin{eqnarray*}
\left|\sum_{i=1}^n A_i\right|^r &=& \left(\left|\sum_{i=1}^n A_i\right|^2\right)^{r/2}\\
&=&\left(\sum_{i,j=1}^n A_i^*A_j\right)^{r/2}\\
&=&\left(\sum_{i=1}^n |A_i|^2\right)^{r/2}\\
&\leq&\left(\sum_{i=1}^n
\alpha_i^{2/(2-r)}\right)^{(r-2)/2}\sum_{i=1}^n
\alpha_i\left(|A_i|^2\right)^{r/2}\qquad
({\rm by~Corollary~ \ref{xi}~with~} X_i=I)\\
&=&\left(\sum_{i=1}^n
\alpha_i^{2/(2-r)}\right)^{(r-2)/2}\sum_{i=1}^n
\alpha_i\left|A_i\right|^r\;.
\end{eqnarray*}
Inequality \eqref{norm} is easily deduced from the fact that
$\|Z\|^r=\|\, |Z|^r\|$ for each $Z\in {\mathbb B}({\mathscr H})$.
\end{proof}

\begin{remark}
It is clear that $A_1, \cdots A_n \in {\mathbb B}({\mathscr H})$
have orthogonal ranges if and only if $A_i^*A_j=0$. An example of
such operators is obtained by considering an orthogonal family
$(e_i)_{1\leq i\leq n}$ and a vector $x$ in ${\mathscr H}$ and
defining the rank one operators $A_i: {\mathscr H} \to {\mathscr H}$
by $A_i=e_i \otimes x$, $1 \leq i \leq n$. Then $A_i^*A_j=\langle
e_j, e_i\rangle\, x\otimes x$ for all $1 \leq i,j \leq n$.
\end{remark}

\end{document}